\documentclass{article}
\usepackage{amssymb,bm}
\usepackage{amsmath}
\usepackage{amsthm}
\usepackage{amsfonts}
\usepackage{mathrsfs}
\usepackage{appendix}
\usepackage{hyperref}
\usepackage{authblk}

\newtheorem{thm}{Theorem}[section]
\newtheorem{lem}[thm]{Lemma}

\theoremstyle{remark}

\usepackage{graphics}
\usepackage{graphicx}
\usepackage{geometry}
\usepackage{algorithm}
\usepackage{algorithmic}
\hypersetup{colorlinks,linkcolor={blue},citecolor={blue},urlcolor={red}}

\usepackage[english]{babel}

\usepackage{cite}
\usepackage{threeparttable}
\usepackage{float}
\floatstyle{plaintop}
\restylefloat{table}

\floatstyle{plain}
\restylefloat{figure}
\usepackage{placeins}

\def\x{{\mathbf x}}
\def\y{{\mathbf y}}

\begin{document}

\title{Necessary and Sufficient Null Space Condition for Nuclear Norm Minimization in Low-Rank Matrix Recovery}

\author{
Jirong Yi \textsuperscript{a}{\footnote{\textsuperscript{a}
Department of Electrical and Computer Engineering, University of Iowa, Iowa City, IA 52242.}}\,
Weiyu Xu\textsuperscript{c}\thanks{\noindent\textsuperscript{c}Department of Electrical and Computer Engineering, University of Iowa, Iowa City, IA 52242. \\
Email: \texttt{weiyu-xu@uiowa.edu}}
}

\maketitle

\begin{abstract}

Low-rank matrix recovery has found many applications in science and engineering such as machine learning, signal processing, collaborative filtering, system identification, and Euclidean embedding. But the low-rank matrix recovery problem is an NP hard problem and thus challenging. A commonly used heuristic approach is the nuclear norm minimization. In \cite{recht_null_2011,recht_necessary_2008,oymak_new_2010}, the authors established the necessary and sufficient null space conditions for nuclear norm minimization to recover every possible low-rank matrix with rank at most $r$ (the strong null space condition). In addition, in \cite{oymak_new_2010}, Oymak et al. established a null space condition for successful recovery of a given low-rank matrix (the weak null space condition) using nuclear norm minimization, and derived the phase transition for the nuclear norm minimization. In this paper, we show that the weak null space condition in \cite{oymak_new_2010} is only a sufficient condition for successful matrix recovery using nuclear norm minimization, and is not a necessary condition as claimed in \cite{oymak_new_2010}. In this paper, we further give a weak null space condition for low-rank matrix recovery, which is both necessary and sufficient for the success of nuclear norm minimization. At the core of our derivation are an inequality for characterizing the nuclear norms of block matrices, and the conditions for equality to hold in that inequality.
\end{abstract}


\section{INTRODUCTION}

The problem of reconstructing matrices from limited observations has attracted significant attention in many areas such as machine learning \cite{candes_exact_2009,candes_power_2010,cai_singular_2010}, computer vision \cite{tomasi_shape_1992}, and phaseless signal recovery \cite{candes_phase_2013}. Very often we deal with the problem of finding a low-rank matrix consistent with existing observations, known as the rank minimization problem (RMP) \cite{recht_guaranteed_2010,fazel_matrix_2002}.

Let $X\in\mathbb{R}^{n_1\times n_2}$ be the ground truth low-rank matrix, and we observe $X$ through a linear mapping $\mathcal{A}: \mathbb{R}^{n_1\times n_2} \rightarrow \mathbb{R}^m$. Suppose the measurement result is $b=\mathcal{A}(X)$, then the RMP can be formulated as follows:
\begin{align}\label{Defn:RankMinimizationProblemIntr}
&{\rm minimize\ } {\rm rank}(Y)\\
&{\rm subject\ to\ }\mathcal{A}(Y)=b.
\end{align}
Because the RMP is an NP-hard problem, researchers often relax it to the nuclear norm minimization (NNM):
\begin{align}\label{Defn:NuclearNormMinimizationIntr}
&{\rm minimize\ } \|Y\|_*\nonumber\\
&{\rm subject\ to\ }\mathcal{A}(Y)=b,
\end{align}
where $\|\cdot\|_*$ is the nuclear norm, namely the sum of the singular values of a matrix. Over the years, there have been a large volume of research results on nuclear norm minimization (\ref{Defn:NuclearNormMinimization}), including deriving recovery performance guarantees \cite{recht_guaranteed_2010,candes_exact_2009,candes_power_2010} and designing numerical methods for solving it \cite{liu_interior-point_2009,ma_fixed_2011,cai_singular_2010,li_reweighted_2014}.

The properties of $\mathcal{A}$ play an important role in establishing the recovery guarantees of nuclear norm minimization. The restricted isometry property (RIP) for $\mathcal{A}$ is often used to prove performance guarantees of the nuclear norm minimization \cite{recht_guaranteed_2010, candes_exact_2009}. However, the RIP is only a suffcient but not a necessary condition for the success of nuclear norm minimization \cite{recht_null_2011}.

In \cite{recht_null_2011,recht_necessary_2008}, the authors characterized the necessary and sufficient null space condition for successful reconstruction of every ground truth matrix with rank no more than $r$ (the strong null space condition). In \cite{oymak_new_2010}, Oymak et al. obtained a more concise and verifiable necessary and sufficient null space condition in the strong sense. The authors established a null space condition for successful recovery of a given low-rank matrix (the weak null space condition) in \cite{oymak_new_2010}. Compared with the strong null space condition, the weak null space condition is a condition on the null space of $\mathcal{A}$ such that a particular ground truth rank-$r$ matrix can be successfully recovered using nuclear norm minimization. In \cite{oymak_new_2010}, it was claimed that the weak null space condition discovered in \cite{oymak_new_2010}  was both a necessary and sufficient condition for recovering a particular rank-$r$ matrix using nuclear norm minimization.

However, in this paper, we show that the weak null space condition proposed in \cite{oymak_new_2010} is not a necessary condition for the success of nuclear norm minimization in low-rank matrix recovery. Furthermore, we provide a new derivation which gives a true necessary and sufficient weak null space condition for the nuclear norm minimization. At the core of our derivation are an inequality for characterizing the nuclear norms of block matrices, and the conditions for equality to hold in that inequality.

This paper is organized as follows. In Section \ref{Sec:PreliminariesAndProblemFormulation}, we formulate the matrix recovery problem as a nuclear norm minimization problem. We provide a counterexample to illustrate that the weak null space condition in \cite{oymak_new_2010} is not a necessary condition. In Section \ref{Sec:NullSpaceCondition}, we give formal statements of our main theorems and their proofs. We also give key lemmas needed for proving the main results. We present the proofs of these lemmas in Appendix.

\textbf{Notations}: In this paper, we denote the nuclear norm of a matrix $X$ by $\|X\|_*=\sum_{i=1} \sigma_i(X)$, where
$\sigma_i(X)$ is the $i$-th largest singular value of $X$. The bold $\bm{0}$ is used to denote all-zero vectors or all-zero matrices, and its dimension depends on the context. The trace of a matrix $X$ is denoted by ${\rm Tr}(X)$, and the transpose and conjugate transpose of a matrix $X$ are denoted by $X^T$ and $X^*$ respectively.

\section{PRELIMINARIES}\label{Sec:PreliminariesAndProblemFormulation}

Let $X\in\mathbb{R}^{n_1\times n_2} (n_1\leq n_2)$ be of rank $r\leq\min\{n_1,n_2\}$, $\mathcal{A}: \mathbb{R}^{n_1\times n_2}\rightarrow \mathbb{R}^m$ be a linear measruement operator, and $b=\mathcal{A}(X)\in\mathbb{R}^m$ be the measurement vector. We consider the following rank minimization problem:
\begin{align}\label{Defn:RankMinimizationProblem}
&{\rm minimize\ }{\rm rank}(Y)\nonumber\\
&{\rm subject\ to\ }\mathcal{A}(Y)=b.
\end{align}

Recently, in \cite{oymak_new_2010}, Oymak et al. proposed the following weak null space condition as a ``necessary'' and sufficient condition for the nuclear norm minimization (\ref{Defn:NuclearNormMinimization}):
\begin{align}\label{Defn:NuclearNormMinimization}
&{\rm minimize\ }\|Y\|_*\nonumber\\
&{\rm subject\ to\ }\mathcal{A}(Y)=b.
\end{align}

{\bf Weak null space condition from \cite{oymak_new_2010}:} Let $X\in\mathbb{R}^{n\times n}$ be a matrix with rank $r$, and let its singular value decomposition (SVD) be $X=U\Sigma V^T$ with $\Sigma\in\mathbb{R}^{r\times r}$. Then $X$ will be the unique solution to  (\ref{Defn:NuclearNormMinimization}) if and only if for all nonzero $W\in\mathcal{N}(\mathcal{A})$, we have
\begin{align}\label{WeakNullSpaceCondition_Oymak}
{\rm Tr}(U^TWV)+\|\bar{U}^TW\bar{V}\|_*>0,
\end{align}
where $\bar{U}$ and $\bar{V}$ are chosen such that $[U\ \bar{U}]$ and $[V\ \bar{V}]$ are unitary, and $\mathcal{N}(\mathcal{A})$ is the null space of $\mathcal{A}$.

 However, we discover that the above weak null space condition from \cite{oymak_new_2010} is only a sufficient but not necessary condition for the success of nuclear norm minimization. Here we present a simple counterexample where (\ref{WeakNullSpaceCondition_Oymak}) is violated, but the nuclear norm minimization still succeeds in recovering the ground truth matrix $X$. To simplify presentation, we first use a counterexample in the field of real numbers (where every element in the null space is a real-numbered matrix
  ) to illustrate the idea. Building on this real-numbered example, we further give a counterexample in the field of complex numbers.

Suppose \begin{align}
{X}=
\begin{array}{l}
\left[\begin{array}{*{20}{c}}
-1 &0\\
0 &0
\end{array}\right]
\end{array},~
{Q}=
\begin{array}{l}
\left[\begin{array}{*{20}{c}}
1 &1\\
1 &1
\end{array}\right].
\end{array}
\end{align}
We also assume that the linear mapping $\mathcal{A}$ is such that ${tQ} (t\neq0)$ is the only type of nonzero elements in the null space of $\mathcal{A}$. Then the solution to (\ref{Defn:NuclearNormMinimization}) must be of the form ${X}+t{Q}$, where $t$ is any real number. Let the singular value decomposition of $X$ be $X=U\Lambda V^*$, where
\begin{align}
{U}=
\begin{array}{l}
\left[\begin{array}{*{20}{c}}
1\\ 0
\end{array}\right]
\end{array}
,
{V}=
\begin{array}{l}
\left[\begin{array}{*{20}{c}}
-1\\ 0
\end{array}\right]
\end{array}
.
\end{align}
Define $\bar{U}$ and $\bar{V}$ as
\begin{align}
\bar{{U}}=
\begin{array}{l}
\left[\begin{array}{*{20}{c}}
0\\ 1
\end{array}\right]
\end{array}
,
\bar{{V}}=
\begin{array}{l}
\left[\begin{array}{*{20}{c}}
0\\ 1
\end{array}\right]
\end{array}.
\end{align}

One can check that, for this example,
\begin{align}
-|{\rm Tr}({U}^*{Q}{V})|+\|\bar{{U}}^*{Q}\bar{{V}}\|_*= 0.
\end{align}
However, we will show that
\begin{align}
\left\vert \left\vert {X}+t{Q} \right\vert\right\vert_*>1, \forall t\neq0,
\end{align}
implying that ${X}$ is the unique solution to (\ref{Defn:NuclearNormMinimization}).

In fact, we calculate
\begin{align}
{B}
=({X}+t{Q})({X}+t{Q})^T
=
\begin{array}{l}
\left[\begin{array}{*{20}{c}}
(-1+t)^2+t^2 & (-1+t)t+t^2\\
(-1+t)t+t^2 & 2t^2
\end{array}\right]
\end{array},
\end{align}
and then the singular values of ${X}+t{Q}$ are the square roots of the eigenvalues of ${B}$. The eigenvalues of ${B}$ can be obtained by solving for $\lambda$ using
\begin{align}
{\rm det}({B}-\lambda {I})=0,
\end{align}
where $I$ is the identity matrix and ${\rm det}(\cdot)$ is the determinant. This results in
\begin{align}
\lambda=\frac{a(t)+b(t)\pm\sqrt{(a(t)-b(t))^2+4c(t)}}{2},
\end{align}
where
\begin{align}
a(t)=(-1+t)^2+t^2,
b(t)=2t^2,
c(t)=((-1+t)t+t^2)^2.
\end{align}
Thus the two eigenvalues of ${B}$ are
\begin{align}
&\lambda_1=\frac{4t^2-2t+1+\sqrt{16t^4-16t^3+8t^2-4t+1}}{2},\\
&\lambda_2=\frac{4t^2-2t+1-\sqrt{16t^4-16t^3+8t^2-4t+1}}{2},
\end{align}
and the singular values of ${X}+t{Q}$ are
\begin{align}
&\sigma_1=\sqrt{\lambda_1}=
\sqrt{\frac{4t^2-2t+1+\sqrt{16t^4-16t^3+8t^2-4t+1}}{2}},\\
&\sigma_2=\sqrt{\lambda_2}=
\sqrt{\frac{4t^2-2t+1-\sqrt{16t^4-16t^3+8t^2-4t+1}}{2}}.
\end{align}
After some algebra, we get
\begin{align}
\|{X}+t{Q}\|_*
=\sigma_1+\sigma_2
&=
\begin{cases}
\sqrt{4t^2+1}, t\geq0,\\
1-2t, t<0.
\end{cases}
\end{align}
This means $\|{X}+t{Q}\|_*$  is always greater than $1$ for $t\neq0$, showing ${X}$ is the unique solution to the nuclear norm minimization. But $-|{\rm Tr}({U}^*{Q}{V})|+\|\bar{{U}}^*{Q}\bar{{V}}\|_*\ngtr 0.$

Now we give a counterexample in the field of complex numbers, where the null space of $\mathcal{A}$ contains complex-numbered matrices. Suppose that we have the same
matrices ${X}$ and ${Q}$. Then the solution to (\ref{Defn:NuclearNormMinimization}) must be of the form ${X}+t{Q}$, where $t$ is any complex number. Without loss of generality, let us
take $t=-ae^{-\imath\theta}$, where $a\geq 0$ is a nonnegative real number, $\theta$ is any real number between 0 and $2\pi$, and $\imath = \sqrt{-1}$. We further denote $
{B}
=({X}+t{Q})({X}+t{Q})^*$. Then by calculating the eigenvalues of ${B}$, we obtain that
\begin{align}
&\|{X}+t{Q}\|_*=\sigma_1+\sigma_2=\sqrt{4a^2+2a(1+\cos(\theta))+1}\\
&=\sqrt{4 \left (a+\frac{1+\cos(\theta)}{4}\right)^2+1-\frac{(1+\cos(\theta))^2}{4}}.
\end{align}
So $\|{X}+t{Q}\|_*>1$, if $a\neq0$ (namely $t\neq0$), implying that the nuclear norm minimization
can uniquely recovers ${X}$ even though  $-|{\rm Tr}({U}^*{Q}{V})|+\|\bar{{U}}^*{Q}\bar{{V}}\|_*\ngtr 0.$

Our counterexample raises the following question: what is a necessary and sufficient weak null space condition for the success of nuclear norm minimization? In the next section, we will answer this question.

\section{NECESSARY AND SUFFICIENT NULL SPACE CONDITION}\label{Sec:NullSpaceCondition}

In this section, we give a necessary and sufficient null space condition for successful recovery of the ground truth matrix using nuclear norm minimization. Our main results are stated in Theorems \ref{Thm:SharpWeakNullSpaceCondition} and \ref{Thm:SharpWeakNullSpaceConditionComplexCase}. To simplify our presentations, in Theorem \ref{Thm:SharpWeakNullSpaceCondition} and its proof, we only consider the case of $X$ being a square real-numbered matrix. In Theorem \ref{Thm:SharpWeakNullSpaceConditionComplexCase}, we generalize our results to complex-numbered non-square matrices without proof.

\begin{thm}\label{Thm:SharpWeakNullSpaceCondition}
Consider a rank-$r$ matrix $X\in\mathbb{R}^{n\times n}$, and let its singular value decomposition be $X=U_X\Lambda_XV_X^*$. Suppose that we observe $b=\mathcal{A}(X)$, and use the nuclear norm minimization to recover the matrix. For any matrix $Q\in\mathbb{R}^{n\times n}$, we define a matrix $Q'=
\begin{array}{l}
\left[\begin{array}{*{20}{c}}
A' &B'\\
C' &D'
\end{array}\right]
\end{array}\in\mathbb{R}^{n\times n}$ satisfying
\begin{align}
Q=
[U_X\ \overline{U_X}]\times Q'\times [V_X\ \overline{V_X}]^*,
\end{align}
where $[U_X\ \overline{U_X}]$ and $[V_X\ \overline{V_X}]$ are unitary matrices in $\mathbb{R}^{n\times n}$, $A'\in\mathbb{R}^{r\times r}$ and $D'\in\mathbb{R}^{(n-r)\times(n-r)}$. We let the singular value decomposition of $D'$ be $D'=U_{D'}\Lambda_{D'}V_{D'}^*$. We further define a matrix $Q''
=
\begin{array}{l}
\left[\begin{array}{*{20}{c}}
A'' &B''\\
C'' &D''
\end{array}\right]
\end{array}
\in\mathbb{R}^{n\times n}$ such that
\begin{align}
Q'=
\begin{array}{l}
\left[\begin{array}{*{20}{c}}
I_{r\times r} &\bm{0} &\bm{0}\\
\bm{0} & U_{D'} &\overline{U_{D'}}
\end{array}\right]
\end{array}
Q''
\begin{array}{l}
\left[\begin{array}{*{20}{c}}
I_{r\times r} &\bm{0} &\bm{0}\\
\bm{0} &V_{D'} &\overline{V_{D'}}
\end{array}\right]
\end{array}^*,
\end{align}
where $[U_{D'}\ \overline{U_{D'}}]$ and $[V_{D'}\ \overline{V_{D'}}]$ are unitary matrices in $\mathbb{R}^{(n-r)\times (n-r)}, A''\in\mathbb{R}^{r\times r}$ and $D''\in\mathbb{R}^{(n-r)\times (n-r)}$.

Then $X$ is the unique solution to nuclear norm minimization if and only if every nonzero element $Q\in\mathbb{R}^{n\times n}$ from the null space of $\mathcal{A}$ satisfies one of the following two conditions:\\
(1)
\begin{align}
{\rm Tr}(U_X^TQV_X)+\|\overline{U_X}^TQ\overline{V_X}\|_*=0.
\end{align}
Moreover, $Q''$ satisfies at least one of the following conditions: a) the row space of $B''$ is not a subspace of the row space of $D''$; b) the column space of $C''$ is not a subspace of the column space of $D''$; c) $Q''$ is not a symmetric matrix.\\
(2)
\begin{align}
{\rm Tr}(U_X^TQV_X)+\|\overline{U_X}^TQ\overline{V_X}\|_*>0.
\end{align}
\end{thm}

Our proof for Theorems \ref{Thm:SharpWeakNullSpaceCondition} and \ref{Thm:SharpWeakNullSpaceConditionComplexCase} depends on the characterization of the subdifferential for the nuclear norm of real-numbered and complex-numbered matrices \cite{xu_separation-free_2017,cai_robust_2016,usevich_hankel_2016}. Next, we give Lemmas \ref{Lem:NuclearNormEqualityCharacterization-SingularSpace}, \ref{lem:NuclearNormEqualityCharacterization-RowAndColumnSpace}, and \ref{Lem:NulcearNormAndTrace} which play an important role in establishing Theorem \ref{Thm:SharpWeakNullSpaceCondition}.

\begin{lem}\label{Lem:NuclearNormEqualityCharacterization-SingularSpace}
Let us assume that $
X=
\begin{array}{l}
\left[\begin{array}{*{20}{c}}
A &B\\
C &D
\end{array}\right]
\end{array}
$ is a square matrix in $\mathbb{R}^{n\times n}$, where
$A\in \mathbb{R}^{m\times m}$ and $D\in\mathbb{R}^{(n-m)\times (n-m)}$. Let $A=U_A\Lambda_A V_A^*$ and $D=U_D\Lambda_D V_D^*$ be singular value decompositions of $A$ and $D$ respectively. Let $X=U_X\Lambda_XV_X^*$ be the singular value decomposition of $X$.

Then it always holds that
\begin{align}\label{Eq:Inequality}
\left\|
\begin{array}{l}
\left[\begin{array}{*{20}{c}}
A & B\\
C &D
\end{array}\right]
\end{array}
\right\|_*
\geq
\left\|
\begin{array}{l}
\left[\begin{array}{*{20}{c}}
A &\bm{0}\\
\bm{0} &D
\end{array}\right]
\end{array}
\right\|_*.
\end{align}
Moreover, the equality holds if and only if there exists a matrix $Q$ such that $Q^*Q=I$,
\begin{align}\label{Eq:EqualityCondition}
U_X=
\begin{array}{l}
\left[\begin{array}{*{20}{c}}
U_A &\bm{0}\\
\bm{0} &U_D
\end{array}\right]
\end{array}
Q,
\end{align}
and
\begin{align}
V_X=
\begin{array}{l}
\left[\begin{array}{*{20}{c}}
V_A &\bm{0}\\
\bm{0} &V_D
\end{array}\right]
\end{array}
Q.
\end{align}
\end{lem}

We remark that in \cite{recht_null_2011,recht_necessary_2008}, the authors established (\ref{Eq:Inequality}) without specifying the conditions under which (\ref{Eq:Inequality}) takes strict inequality or equality.

\begin{lem}\label{lem:NuclearNormEqualityCharacterization-RowAndColumnSpace}
Let $
X=
\begin{array}{l}
\left[\begin{array}{*{20}{c}}
A &B\\
C &D
\end{array}\right]
\end{array}$ be a square matrix in $\mathbb{R}^{n\times n}
$, where $A\in \mathbb{R}^{m\times m}$ and $D\in\mathbb{R}^{(n-m)\times (n-m)}$. Let $X=U_X\Lambda_XV_X^*$ be the singular value decomposition of $X$, and let $A=U_A\Lambda_A V_A^*$ and $D=U_D\Lambda_D V_D^*$ be the singular value decomposition of $A$ and $D$. Let $
[U_A\ \overline{U_A}]$ and
$[V_A\ \overline{V_A}]$ be unitary matrices in $\mathbb{R}^{m\times m}$, and let $
[U_D\ \overline{U_D}]$ and $
[V_D\ \overline{V_D}]$ be unitary matrices in $\mathbb{R}^{(n-m)\times (n-m)}$. Let us define
\begin{align}
E=
\begin{array}{l}
\left[\begin{array}{*{20}{c}}
U_A &\overline{U_A} &\bm{0} &\bm{0}\\
\bm{0} &\bm{0} &U_D &\overline{U_D}
\end{array}\right]
\end{array}^*
\times
\begin{array}{l}
\left[\begin{array}{*{20}{c}}
A &B\\
C &D
\end{array}\right]
\end{array}
\times
\begin{array}{l}
\left[\begin{array}{*{20}{c}}
V_A &\overline{V_A} &\bm{0} &\bm{0}\\
\bm{0} &\bm{0} &V_D &\overline{V_D}
\end{array}\right].
\end{array}
\end{align}
Then
\begin{align}
\left\|
\begin{array}{l}
\left[\begin{array}{*{20}{c}}
A &B\\
C &D
\end{array}\right]
\end{array}
\right\|_*
= \left\|
\begin{array}{l}
\left[\begin{array}{*{20}{c}}
A &\bm{0}\\
\bm{0} &D
\end{array}\right]
\end{array}
\right\|_*
\end{align}
if and only if the following conditions hold simultaneously: (1) the row space of $X$ is a subspace of the row space of $\begin{array}{l}
\left[\begin{array}{*{20}{c}}
V_A &\bm{0}\\
\bm{0} &V_D
\end{array}\right]
\end{array}^*$; (2) the column space of $X$ is a subspace of the column space of $\begin{array}{l}
\left[\begin{array}{*{20}{c}}
U_A &\bm{0}\\
\bm{0} &U_D
\end{array}\right]
\end{array}$; (3) $E$ is a symmetric positive semidefinite matrix.
\end{lem}

\begin{lem}\label{Lem:NulcearNormAndTrace}
Let $A$ be an $n\times n$ matrix. Then
\begin{align}
\|A\|_*\geq\sum_{i=1}^n |A_{i,i}|.
\end{align}
Moreover,
\begin{align}
\|A\|_*=\sum_{i=1}^n A_{i,i}
\end{align}
if and only if $A$ is a symmstric positive semidefinite matrix.
\end{lem}
We note that Lemma \ref{Lem:NulcearNormAndTrace} improves over Lemma 3 in \cite{oymak_new_2010}. With the lemmas above, we now present the proof for Theorem \ref{Thm:SharpWeakNullSpaceCondition}.

\emph{Proof for Theorem \ref{Thm:SharpWeakNullSpaceCondition}:}
We define a function $f(X)=\|X\|_*$. We look at three cases for the sign of ${\rm Tr}(U_X^TQV_X)+\|\overline{U_X}^TQ\overline{V_X}\|_*$: a) negative, b) positive, and c) being equal to 0.
~\\

\noindent\textbf{Case a)}:\\
We first assume that, for a certain nonzero $Q$ from the null space of $\mathcal{A}$,
\begin{align}\label{Eq:contrarysmaller}
{\rm Tr}(U_X^TQV_X)+\|\overline{U_X}^TQ\overline{V_X}\|_*<0.
\end{align}
Because
\begin{align}\label{Eq:contrary1}
{\rm Tr}(U_X^TQV_X)+\|\overline{U_X}^TQ\overline{V_X}\|_*
=\sup_{W \in \partial \|X\|_*} <W, Q> ,
\end{align}
we have
$\sup_{W \in \partial \|X\|_*} <W, Q>\ <0$.

From Theorem 23.4 in \cite{rockafellar_convex_2015}, we know that
\begin{align}\label{Eq:contrary2}
f'(X;Q)=\sup_{W \in \partial \|X\|_*} <W, Q>\ <0,
\end{align}
where
\begin{align}\label{Eq:contrary3}
f'(X;Q)=\inf_{t>0} \frac{\|X+tQ\|_*-\|X\|_*}{t}
\end{align}
is the one-sided directional derivative defined in \cite{rockafellar_convex_2015}.  Then there must exist
a positive number $\lambda$ such that $\|X+\lambda Q\|_*-\|X\|_*<0$, namely $Q$ is a direction along which we can reduce the nuclear norm of $X$. This shows that $X$ cannot be the solution to the nuclear norm minimization problem, if ${\rm Tr}(U_X^TQV_X)+\|\overline{U_X}^TQ\overline{V_X}\|_*<0$ for a certain nonzero $Q$ from the null space of $\mathcal{A}$.

~\\

\noindent\textbf{Case b):}\\
Now we instead assume that, a nonzero matrix $Q$ from the null space of $\mathcal{A}$ satisfies ${\rm Tr}(U_X^TQV_X)+\|\overline{U_X}^TQ\overline{V_X}\|_*>0$. We know that the subdifferential $\partial f(X)$ of $f(X)$ is given by the set of matrices in the form of $U_XV_X^T+ \overline{U_X}M \overline{V_X}^T$, where $M$ is any matrix (of appropriate dimension) with spectral norm no more than $1$. Then by the definition of subdifferential, we have
\begin{align}
\|X+Q\|_*
& \geq \|X\|_*+\sup_{B \in \partial f(X)} \langle Q, B \rangle    \\
&=\|X\|_*+{\rm Tr}(U_X^TQV_X)+\|\overline{U_X}^TQ\overline{V_X}\|_*\\
&> \|X\|_*,
\end{align}
where we use the fact that the dual norm of spectral norm is the nuclear norm. This implies that $X+tQ$ cannot be the solution to nuclear norm minimization for $t\neq 0$ and any nonzero $Q$ from the null space of $\mathcal{A}$ satisfying ${\rm Tr}(U_X^TQV_X)+\|\overline{U_X}^TQ\overline{V_X}\|_*>0$.

~\\
\noindent\textbf{Case c:)}\\
Now, we only need to consider the case ${\rm Tr}(U_X^*QV_X) + \|\overline{U_X}^*Q\overline{V_X}\|_*=0$. We know that any matrix $Y$ satisfying $b=\mathcal{A}(X)$ must be of the form
\begin{align}
Y=X+tQ,
\end{align}
where $t>0$, and $Q$ is a nonzero element in the null space of $\mathcal{A}$.

Let us consider a matrix $Q$ in the null space of $\mathcal{A}$ such that
\begin{align}
{\rm Tr}(U^TQV)+\|\bar{U}^TQ\bar{V}\|_*=0.
\end{align}

Let $[U_X\ \overline{U_X}]$ and $[V_X\ \overline{V_X}]$ be unitary matrices in $\mathbb{R}^{n\times n}$. Then we can express $Q$ as
\begin{align}
Q=[U_X\ \overline{U_X}]\times Q'\times [V_X\ \overline{V_X}]^*,
\end{align}
where $Q'\in\mathbb{R}^{n\times n}$.

We write $Q'$ as a block matrix
\begin{align}
Q'=
\begin{array}{l}
\left[\begin{array}{*{20}{c}}
A' &B'\\
C' &D'
\end{array}\right]
\end{array},
\end{align}
where $A'\in\mathbb{R}^{r\times r}$ and $D'\in\mathbb{R}^{(n-r)\times (n-r)}$.

We let the singular value decomposition of $D'$ be $D'=U_{D'}\Lambda_{D'}V_{D'}^*$. Then we further express $Q'$ as
\begin{align}
Q'=
\begin{array}{l}
\left[\begin{array}{*{20}{c}}
I_{r\times r} &\bm{0}\\
\bm{0} &
{\begin{array}{l}
\begin{array}{*{20}{c}}
U_{D'} &\overline{U_{D'}}
\end{array}
\end{array}}
\end{array}\right]
\end{array}
\times Q'' \times
\begin{array}{l}
\left[\begin{array}{*{20}{c}}
I_{r\times r} &\bm{0}\\
\bm{0} &
{\begin{array}{l}
\begin{array}{*{20}{c}}
V_{D'} &\overline{V_{D'}}
\end{array}
\end{array}}
\end{array}\right]
\end{array}^*
\end{align}
where $Q''\in\mathbb{R}^{n\times n}$, and $[U_{D'}\ \overline{U_{D'}}]$ and $[V_{D'}\ \overline{V_{D'}}]$ are unitary matrices in $\mathbb{R}^{(n-r)\times (n-r)}$.

We can write $Q''$ as a block matrix
\begin{align}
Q''=
\begin{array}{l}
\left[\begin{array}{*{20}{c}}
A'' &B''\\
C'' &D''
\end{array}\right]
\end{array}
\end{align}
where $A''\in\mathbb{R}^{r\times r}, D''\in\mathbb{R}^{(n-r)\times(n-r)}$. Moreover,
\begin{align}
D''=
\begin{array}{l}
\left[\begin{array}{*{20}{c}}
\Lambda_{D'} &\bm{0}\\
\bm{0} &\bm{0}
\end{array}\right]
\end{array}
\in\mathbb{R}^{(n-r)\times (n-r)}.
\end{align}

From the singular value decomposition of $X$, we also have
\begin{align}
X& =
\begin{array}{l}
\left[\begin{array}{*{20}{c}}
U_X &\overline{U_X}
\end{array}\right]
\end{array}
\times
\begin{array}{l}
\left[\begin{array}{*{20}{c}}
I_{r\times r} &\bm{0}\\
\bm{0} &
{\begin{array}{l}
\begin{array}{*{20}{c}}
U_{D'} &\overline{U_{D'}}
\end{array}
\end{array}}
\end{array}\right]
\end{array}
\times
\begin{array}{l}
\left[\begin{array}{*{20}{c}}
\Lambda_X &\bm{0}\\
\bm{0} &\bm{0}
\end{array}\right]
\end{array}\nonumber\\
&\ \ \ \ \times \begin{array}{l}
\left[\begin{array}{*{20}{c}}
I_{r\times r} &\bm{0}\\
\bm{0} &
{\begin{array}{l}
\begin{array}{*{20}{c}}
V_{D'} &\overline{V_{D'}}
\end{array}
\end{array}}
\end{array}\right]
\end{array}^*
\times
\begin{array}{l}
\left[\begin{array}{*{20}{c}}
V_X & \overline{V_X}
\end{array}\right]
\end{array}^*.
\end{align}
Then
\begin{align}
Y=X+tQ
& =
\begin{array}{l}
\left[\begin{array}{*{20}{c}}
U_X &\overline{U_X}
\end{array}\right]
\end{array}
\times
\begin{array}{l}
\left[\begin{array}{*{20}{c}}
I_{r\times r} &\bm{0}\\
\bm{0} &
{\begin{array}{l}
\begin{array}{*{20}{c}}
U_{D'} &\overline{U_{D'}}
\end{array}
\end{array}}
\end{array}\right]
\end{array}
\times Y''\nonumber\\
& \ \ \ \ \times
\begin{array}{l}
\left[\begin{array}{*{20}{c}}
I_{r\times r} &\bm{0}\\
\bm{0} &
{\begin{array}{l}
\begin{array}{*{20}{c}}
V_{D'} &\overline{V_{D'}}
\end{array}
\end{array}}
\end{array}\right]
\end{array}^*
\times
\begin{array}{l}
\left[\begin{array}{*{20}{c}}
V_X & \overline{V_X}
\end{array}\right]
\end{array}^*
\end{align}
where
\begin{align}
Y''=
\begin{array}{l}
\left[\begin{array}{*{20}{c}}
\Lambda_X &\bm{0}\\
\bm{0} &\bm{0}
\end{array}\right]
\end{array}
+
t
\begin{array}{l}
\left[\begin{array}{*{20}{c}}
A'' &B''\\
C'' &D''
\end{array}\right]
\end{array}.
\end{align}

This means that $Y''$ and $Y$ have the same nuclear norm. So in order for us to see whether there exists a $t>0$ such that $\|Y\|_*=\|X\|_*$, we only need to see whether $\|Y''\|_*=\|X\|_*$.

We observe that
\begin{align}
\|Y''\|_*=\left\| \begin{array}{l}
\left[\begin{array}{*{20}{c}}
\Lambda_X &\bm{0}\\
\bm{0} &\bm{0}
\end{array}\right]
\end{array}
+
t
\begin{array}{l}
\left[\begin{array}{*{20}{c}}
A'' &B''\\
C'' &D''
\end{array}\right]
\end{array}
 \right\|_*.
 \end{align}

Furthermore, because ${\rm Tr}(U^TQV)+\|\overline{U}^TQ\overline{V}\|_*={0}$, we have
\begin{align}
\sum_{i=1}^r A_{i,i}''+\|D''\|_*={0}.
\end{align}

According to Lemma \ref{Lem:NulcearNormAndTrace}, we first note that $\|\Lambda_X+tA''\|_*\geq \sum_{i=1}^r |(\Lambda_X)_{i,i}+tA''_{i,i}|\geq \sum_{i=1}^r ((\Lambda_X)_{i,i}+tA''_{i,i})$. Moreover, by Lemma \ref{Lem:NuclearNormEqualityCharacterization-SingularSpace}, we have
\begin{align*}
\left\|
\begin{array}{l}
\left[\begin{array}{*{20}{c}}
\Lambda_X &\bm{0}\\
\bm{0} &\bm{0}
\end{array}\right]
\end{array}
+t
\begin{array}{l}
\left[\begin{array}{*{20}{c}}
A'' &B''\\
C'' &D''
\end{array}\right]
\end{array}
\right\|_*
& \geq \|\Lambda_X+tA''\|_*+t\|D''\|_*\\
& \geq \|\Lambda_X\|_*+t\sum_{i=1}^r A_{i,i}+t\|D''\|_*\\
& = \|\Lambda_X\|_*\\
& = \|X\|_*.
\end{align*}
Thus if for a certain $t>0$, $\|Y''\|_*=\|X\|_*$, we must have
\begin{align}\label{Eq:SearchForm}
\left\|
\begin{array}{l}
\left[\begin{array}{*{20}{c}}
\Lambda_X &\bm{0}\\
\bm{0} &\bm{0}
\end{array}\right]
\end{array}
+t
\begin{array}{l}
\left[\begin{array}{*{20}{c}}
A'' &B''\\
C'' &D''
\end{array}\right]
\end{array}
\right\|_*
& =
\left\|
\begin{array}{l}
\left[\begin{array}{*{20}{c}}
\Lambda_X+tA'' &\bm{0}\\
\bm{0} &tD''
\end{array}\right]
\end{array}
\right\|_* \nonumber\\
& =
\left\|
\begin{array}{l}
\left[\begin{array}{*{20}{c}}
\Lambda_X+t{\rm diag}(A'') &\bm{0}\\
\bm{0} & tD''
\end{array}\right]
\end{array}
\right\|_*,
\end{align}
where ${\rm diag}(A'')$ is a diagonal matrix having same diagonal as $A''$, and $\Lambda_X+t{\rm diag}(A'')$ has nonnegative diagonal elements.

By Lemma \ref{Lem:NulcearNormAndTrace}, (\ref{Eq:SearchForm}) holds and $\Lambda_X+t{\rm diag}(A'')$ has nonnegative diagonal elements, if and only if
\begin{align}
\begin{array}{l}
\left[\begin{array}{*{20}{c}}
\Lambda_X &\bm{0}\\
\bm{0} &\bm{0}
\end{array}\right]
\end{array}
+t
\begin{array}{l}
\left[\begin{array}{*{20}{c}}
A'' &B''\\
C'' &D''
\end{array}\right]
\end{array}
\end{align}
is a symmetric positive semidefinite matrix for some $t>0$.

We now prove that, there exists $t>0$ such that
\begin{align}
\begin{array}{l}
\left[\begin{array}{*{20}{c}}
\Lambda_X &\bm{0}\\
\bm{0} &\bm{0}
\end{array}\right]
\end{array}
+t
\begin{array}{l}
\left[\begin{array}{*{20}{c}}
A'' &B''\\
C'' &D''
\end{array}\right]
\end{array}
\end{align}
 is a symmetric positive semidefinite matrix if and only if 1) the row space of $B''$ is a subspace of the row space of $D''$; 2) the column space of $C''$ is a subspace of the column space of $D''$; 3) $\begin{array}{l}
\left[\begin{array}{*{20}{c}}
A'' &B''\\
C'' &D''
\end{array}\right]
\end{array}$
is a symmetric matrix.

To see this, we first assume these three conditions are satisfied. Then when $t>0$ is small enough, $
\begin{array}{l}
\left[\begin{array}{*{20}{c}}
\Lambda_X &\bm{0}\\
\bm{0} &\bm{0}
\end{array}\right]
\end{array}
+t\begin{array}{l}
\left[\begin{array}{*{20}{c}}
A'' &B''\\
C'' &D''
\end{array}\right]
\end{array}$ must be a symmetric positive semidefinite matrix. This is because $tD''=t
\begin{array}{l}
\left[\begin{array}{*{20}{c}}
\Lambda_{D'} &\bm{0}\\
\bm{0} & \bm{0}
\end{array}\right]
\end{array}$ is a positive semidefinite matrix, the row space of $tB''$ is a subspace of the row space of $tD''$, the column space of $tC''$ is a subspace of column space of $tD''$, and the Schur complement $\Lambda_X+tA''-(tB'')(t{D}'')^{\dagger}(tC'')$ must be positive semidefinite when $t$ is sufficiently small (noting that $\Lambda_X$ is a full-rank positive semidefinite matrix), where $(tD'')^\dagger$ is the Moore-Penrose inverse of $(tD'')$.

Now we show that if either one of the three conditions fails, then for every $t>0$,
\begin{align}
\begin{array}{l}
\left[\begin{array}{*{20}{c}}
\Lambda_X &\bm{0}\\
\bm{0} &\bm{0}
\end{array}\right]
\end{array}
+t
\begin{array}{l}
\left[\begin{array}{*{20}{c}}
A'' &B''\\
C'' &D''
\end{array}\right]
\end{array}
\end{align}
 is not a symmetric positive semidefinite matrix. This is because, by Schur complement criterion, $\begin{array}{l}\left[\begin{array}{*{20}{c}}
\Lambda_X &\bm{0}\\
\bm{0} &\bm{0}
\end{array}\right]
\end{array}
+t
\begin{array}{l}
\left[\begin{array}{*{20}{c}}
A'' &B''\\
C'' &D''
\end{array}\right]
\end{array}$ is a symmetric positive semidefinite matrix, only if 1) the row space of $tB''$ is a subspace of the row space of $tD''$; 2) the column space of $tC''$ is a subspace of the column space of $tD''$; and 3) $\begin{array}{l}
\left[\begin{array}{*{20}{c}}
A'' &B''\\
C'' &D''
\end{array}\right]
\end{array}$
is a symmetric matrix.

If either of the three conditions fails, we have
\begin{align*}
\left\|
\begin{array}{l}
\left[\begin{array}{*{20}{c}}
\Lambda_X &\bm{0}\\
\bm{0} &\bm{0}
\end{array}\right]
\end{array}
+t
\begin{array}{l}
\left[\begin{array}{*{20}{c}}
A'' &B''\\
C'' &D''
\end{array}\right]
\end{array}
\right\|_*
& >
\|\Lambda_X+tA''\|_*+t\|D''\|_*\\
& \geq \|\Lambda_X\|_*+t\sum_{i=1}^r A_{i,i}+t\|D''\|_*\\
& =\|\Lambda_X\|_* \\
& =\|X\|_*.
\end{align*}
Thus if either of these three conditions fails, the the ground truth $X$ will be the unique solution to the nuclear norm minimization.

Combining Cases a), b) and c), we have proved Theorem \ref{Thm:SharpWeakNullSpaceCondition}.

$\hfill\square$

Now we extend our results to non-square complex-numbered matrices, leading to Theorem \ref{Thm:SharpWeakNullSpaceConditionComplexCase}.  We remark that our proof of Theorem \ref{Thm:SharpWeakNullSpaceConditionComplexCase} depends on the characterization of the subdifferential for the nuclear norm of complex-numbered matrices \cite{xu_separation-free_2017,cai_robust_2016,usevich_hankel_2016}. Other than that, our proof of Theorem \ref{Thm:SharpWeakNullSpaceConditionComplexCase} follows the same line of reasoning as in the proof of Theorem \ref{Thm:SharpWeakNullSpaceCondition}, and so we omit its proof.

\begin{thm}\label{Thm:SharpWeakNullSpaceConditionComplexCase}
Consider a rank-$r$ matrix $X\in\mathbb{C}^{m\times n}$ $(m\leq n)$, and let its singular value decomposition be $X=U_X\Lambda_XV_X^*$. Suppose that we observe $\mathcal{A}(X)$, and use the nuclear norm minimization to recover the matrix. For any matrix $Q\in \mathbb{C}^{m\times n}$, we define a matrix $Q'=
\begin{array}{l}
\left[\begin{array}{*{20}{c}}
A' &B' &\bm{0}\\
C' &D' &\bm{0}
\end{array}\right]
\end{array}\in\mathbb{C}^{m\times n}$ satisfying
\begin{align}
Q=[U_X\ \overline{U_X}]\times Q'\times [V_X\ \overline{V_X}]^*,
\end{align}
where $[U_X\ \overline{U_X}]$ and $[V_X\ \overline{V_X}]$ are unitary matrices in $\mathbb{C}^{m\times m}$ and $\mathbb{C}^{n\times n}$ respectively, and $A'\in\mathbb{C}^{r\times r}$ and
$D'\in\mathbb{C}^{(m-r)\times (m-r)}$. We let the singular value decomposition of $D'$ be $D'=U_{D'}\Lambda_{D'}V_{D'}^*$. We further define $Q''=\begin{array}{l}
\left[\begin{array}{*{20}{c}}
A'' &B'' &\bm{0}\\
C'' &D'' &\bm{0}
\end{array}\right]
\end{array}\in \mathbb{C}^{m\times n}$ such that
\begin{align}
Q'=
\begin{array}{l}
\left[\begin{array}{*{20}{c}}
I_{r\times r} &\bm{0}\\
\bm{0} &
{\begin{array}{l}
\begin{array}{*{20}{c}}
U_{D'} &\overline{U_{D'}}
\end{array}
\end{array}}
\end{array}\right]
\end{array}
\times Q'' \times
\begin{array}{l}
\left[\begin{array}{*{20}{c}}
I_{r\times r} &\bm{0}\\
\bm{0} & {\begin{array}{l}
\begin{array}{*{20}{c}}
V_{D'} &\overline{V_{D'}}
\end{array}
\end{array}}\\
\bm{0} &\bm{0}
\end{array}\right]
\end{array}^*,
\end{align}
where $[U_{D'}\ \overline{U_{D'}}]$ and $[V_{D'}\ \overline{V_{D'}}]$ are unitary matrices in $\mathbb{C}^{(m-r)\times (m-r)}, A''\in\mathbb{C}^{r\times r}$ and $D''\in\mathbb{C}^{(m-r)\times (m-r)}$.

Then $X$ is the unique solution to the nuclea norm minimization if and only if every nonzero element $Q\in\mathbb{C}^{m\times n}$ from the null space of $\mathcal{A}$ satisfies one of the following conditions:\\
1)
\begin{align}
{\rm Re}\{{\rm Tr}(U_X^*QV_X)\}+\|\overline{U_X}^*Q\overline{V_X}\|_*=0.
\end{align}
Moreover, $Q''$ satisfies at least one of the following conditions: a) the row space of $B''$ is not a subspace of the row space of $D''$; b) the column space of $C''$ is not a subspace of the column space of $D''$; c) $\begin{array}{l}
\left[\begin{array}{*{20}{c}}
A'' &B''\\
C'' &D''
\end{array}\right]
\end{array}$ is not a Hermitian matrix.\\
2)
\begin{align}
{\rm Re}\{{\rm Tr}(U_X^*QV_X)\}+\|\overline{U_X}^*Q\overline{V_X}\|_*>0.
\end{align}

\end{thm}

\section{Comparisons with the weak null space condition for sparse recovery}

We would like to contrast the necessary and sufficient weak null space condition for nuclear norm minimization, with the necessary and sufficient weak null space condition for recovery of sparse vectors using $\ell_1$ minimization.  For the $\ell_1$ minimization to successfully (uniquely) recover a sparse vector $\x \in \mathbb{R}^n$ whose support is an index set $K$, the weak null space condition for sparse recovery \cite{donoho_high-dimensional_2006, donoho_neighborliness_2005, xu_precise_2011} requires that for every nonzero element $\y$ in the null space of the linear operator $\mathcal{A}$,
\begin{align}
\label{eq:sparsenull}
\langle\text{sign}(\x_{K}), \y_{K}\rangle +\|\y_{\overline{K}}\|_1>0,
\end{align}
where $\overline{K}=\{1,2,...,n\}\setminus K$.
We remark that one does require strict inequality in the condition (\ref{eq:sparsenull}), in contrast to our new null space condition for nuclear norm minimization where we have shown nuclear norm minimization can still succeed even if equality holds in the null space condition equation in Theorem \ref{Thm:SharpWeakNullSpaceCondition}.   We note that, one can also deduce (\ref{eq:sparsenull}) directly from Theorem \ref{Thm:SharpWeakNullSpaceCondition} by specializing $X$ to be a diagonal matrix with its diagonal elements corresponding to the elements of the vector $\x$. Because none of the three conditions for case ``(1)'' of Theorem \ref{Thm:SharpWeakNullSpaceCondition} holds, one must require strict inequality in case ``2)'' of Theorem  \ref{Thm:SharpWeakNullSpaceCondition}, which is equivalent to (\ref{eq:sparsenull}).

\section*{Appendix}


\section{Proof for Lemma \ref{Lem:NuclearNormEqualityCharacterization-SingularSpace}}

It is known that the dual norm of nuclear norm is the spectral norm. Using this fact, it has been shown in \cite{recht_null_2011} that (\ref{Eq:Inequality}) holds. Now we prove the condition under which the equality holds in (\ref{Eq:Inequality}).

Suppose that the equality in (\ref{Eq:Inequality}) holds. We consider the matrix
\begin{align}
P=
\begin{array}{l}
\left[\begin{array}{*{20}{c}}
U_AV_A^* &\bm{0}\\
\bm{0} &U_DV_D^*
\end{array}\right]
\end{array}.
\end{align}
After some algebra, we know that
\begin{align}
\left\|
\begin{array}{l}
\left[\begin{array}{*{20}{c}}
A &\bm{0}\\
\bm{0} &D
\end{array}\right]
\end{array}
\right\|_*
=
\left<
P,
\begin{array}{l}
\left[\begin{array}{*{20}{c}}
A &{\bm{0}}\\
\bm{0} &D
\end{array}\right]
\end{array}
\right>.
\end{align}
Because we assume that
\begin{align}
\left\|
\begin{array}{l}
\left[\begin{array}{*{20}{c}}
A &B\\
C &D
\end{array}\right]
\end{array}
\right\|_*
=
\left\|
\begin{array}{l}
\left[\begin{array}{*{20}{c}}
A &\bm{0}\\
\bm{0} &D
\end{array}\right]
\end{array}
\right\|_*,
\end{align}
using the fact that the dual norm of nuclear norm is spectral norm, we have
\begin{align}
\sup_{\|Q\|\leq 1} \left<Q,
\begin{array}{l}
\left[\begin{array}{*{20}{c}}
A &B\\
C &D
\end{array}\right]
\end{array}
\right>
=
\left\|
\begin{array}{l}
\left[\begin{array}{*{20}{c}}
A &B\\
C &D
\end{array}\right]
\end{array}
\right\|_*
=
\left<
P,
\begin{array}{l}
\left[\begin{array}{*{20}{c}}
A &\bm{0}\\
\bm{0} &D
\end{array}\right]
\end{array}
\right>
=
\left<
P,
\begin{array}{l}
\left[\begin{array}{*{20}{c}}
A &B\\
C &D
\end{array}\right]
\end{array}
\right>.
\end{align}
Namely $Q=P$ achieves the supremum of
$
\left<Q,
\begin{array}{l}
\left[\begin{array}{*{20}{c}}
A &B\\
C &D
\end{array}\right]
\end{array}
\right>
$ over the set of $Q$'s with spectral norm no more than 1. We have
\begin{align}\label{Eq:ExactSupremum}
\left\|
\begin{array}{l}
\left[\begin{array}{*{20}{c}}
A &B\\
C &D\\
\end{array}\right]
\end{array}
\right\|_*
& =\sup_{\|Q\|\leq 1} \left<Q,
\begin{array}{l}
\left[\begin{array}{*{20}{c}}
A &B\\
C &D
\end{array}\right]
\end{array}
\right>
=\sup_{\|Q\|\leq 1} {\rm Tr}(Q^*U_X\Lambda_X V_X^*) \nonumber\\
& =\sup_{\|Q\|\leq 1} {\rm Tr}(V_X^*Q^*U_X\Lambda_X)
=\sup_{\|Q\|\leq 1} \left<U_X^*QV_X,\Lambda_X\right>\nonumber\\
& =\sup_{\|Q\|\leq 1} \sum_{i=1}^n \sigma_i (U_X^*QV_X)_{ii}
=\sup_{\|Q\|\leq 1} \sum_{i=1}^n \sigma_i U_{X,i:}^* QV_{X,i:}\nonumber\\
& =\sum_{i=1}^n \sigma_i U_{X,i:}^*PV_{X,i:}
\leq \sum_{i=1}^n \sigma_i
\end{align}
where $U_{X,i:}$ and $V_{X,i:}$ are the $i$-th column of $U_X$ and $V_X$ respectively. By the Cauchy-Schwartz inequality, the equality in (\ref{Eq:ExactSupremum}) holds if and only if, for every $i$, the $i$-th column of $U_X$, and the $i$-th column of $V_X$ can be written as $U_{X,i:}=
\begin{array}{l}
\left[\begin{array}{*{20}{c}}
U_A &\bm{0}\\
\bm{0} &U_D
\end{array}\right]
\end{array}
b,
V_{X,i:}=
\begin{array}{l}
\left[\begin{array}{*{20}{c}}
V_A &\bm{0}\\
\bm{0} &V_D
\end{array}\right]
\end{array}
b$ respectively, where $b$ is a vector of unit energy.

Namely, if $\left\|
\begin{array}{l}
\left[\begin{array}{*{20}{c}}
A &B\\
C &D
\end{array}\right]
\end{array}
\right\|_*
=
\left\|
\begin{array}{l}
\left[\begin{array}{*{20}{c}}
A &\bm{0}\\
\bm{0} &D
\end{array}\right]
\end{array}
\right\|_*$, then there exists a matrix $Q$ such that $Q^*Q=I$,
\begin{align}\label{Eq:LeftSingularSpace}
U_X=
\begin{array}{l}
\left[\begin{array}{*{20}{c}}
U_A &\bm{0}\\
\bm{0} &U_D
\end{array}\right]
\end{array}
Q,
\end{align}
and
\begin{align}\label{Eq:RightSingularSpace}
V_X=
\begin{array}{l}
\left[\begin{array}{*{20}{c}}
V_A &\bm{0}\\
\bm{0} &V_D
\end{array}\right]
\end{array}
Q.
\end{align}

On the other hand, if (\ref{Eq:LeftSingularSpace}) and (\ref{Eq:RightSingularSpace}) hold, we have $Q=P=
\begin{array}{l}
\left[\begin{array}{*{20}{c}}
U_AV_A^* &\bm{0}\\
\bm{0} &U_DV_D^*
\end{array}\right]
\end{array}$ achieving the supremum
$$
\sup_{\|Q\|\leq 1}\left<Q,
\begin{array}{l}
\left[\begin{array}{*{20}{c}}
A &B\\
C &D
\end{array}\right]
\end{array}
\right>
$$
through the same arguments in (\ref{Eq:ExactSupremum}). This leads to $\left\|
\begin{array}{l}
\left[\begin{array}{*{20}{c}}
A &B\\
C &D
\end{array}\right]
\end{array}
\right\|_*
=
\left\|
\begin{array}{l}
\left[\begin{array}{*{20}{c}}
A &\bm{0}\\
\bm{0} &D
\end{array}\right]
\end{array}
\right\|_*$.

\section{Proof for Lemma \ref{lem:NuclearNormEqualityCharacterization-RowAndColumnSpace}}

Suppose the rank of $A$ is $r_A$ (possibly smaller than $m$), and the rank of $D$ is $r_D$ (possibly smaller than $(n-m)$). Then we have
\begin{align}
\begin{array}{l}
\left[\begin{array}{*{20}{c}}
A &B\\
C &D
\end{array}\right]
\end{array}
& =
\begin{array}{l}
\left[\begin{array}{*{20}{c}}
U_A\Lambda_AV_A^* &B\\
C &U_D\Lambda_DV_D^*
\end{array}\right]
\end{array} \nonumber\\
& =
\begin{array}{l}
\left[\begin{array}{*{20}{c}}
U_A &\overline{U_A} &\bm{0} &\bm{0}\\
\bm{0} &\bm{0} &U_D &\overline{U_D}
\end{array}\right]
\end{array}
\times
\begin{array}{l}
\left[\begin{array}{*{20}{c}}
{\begin{array}{l}
\begin{array}{*{20}{c}}
\Lambda_A &\bm{0}\\
\bm{0} &\bm{0}
\end{array}
\end{array}} & B'\\
C' &
{\begin{array}{l}
\begin{array}{*{20}{c}}
\Lambda_D &\bm{0}\\
\bm{0} &\bm{0}
\end{array}
\end{array}}
\end{array}\right]
\end{array}\nonumber\\
&\ \ \ \ \times
\begin{array}{l}
\left[\begin{array}{*{20}{c}}
V_A &\overline{V_A} &\bm{0} &\bm{0}\\
\bm{0} &\bm{0} &V_D &\overline{V_D}
\end{array}\right]
\end{array}^*
\end{align}
where $[U_A\ \overline{U_A}]$ and $[V_A\ \overline{V_A}]$ are unitary matrices in $\mathbb{R}^{m\times m}$, $[U_D\ \overline{U_D}]$ and $[V_D\ \overline{V_D}]$ are unitary matrices in $\mathbb{R}^{(n-m)\times (n-m)}$, and $B'$ and $C'$ are such that
\begin{align}
&[U_D \ \overline{U_D}]\ C'\ [V_A \ \overline{V_A}]^*=C,\\
&[U_A \ \overline{U_A}]\ B'\ [V_D \ \overline{V_D}]^*=B.
\end{align}

Since we are performing unitary transformations, we have
\begin{align}
\left\|
\begin{array}{l}
\left[\begin{array}{*{20}{c}}
A &B\\
C &D
\end{array}\right]
\end{array}
\right\|_*
=
\left\|
\begin{array}{l}
\left[\begin{array}{*{20}{c}}
{\begin{array}{l}
\begin{array}{*{20}{c}}
\Lambda_A &\bm{0}\\
\bm{0} &\bm{0}
\end{array}
\end{array}} & B'\\
C' &
{\begin{array}{l}
\begin{array}{*{20}{c}}
\Lambda_D &\bm{0}\\
\bm{0} &\bm{0}
\end{array}
\end{array}}
\end{array}\right]
\end{array}
\right\|_*
=
\|E\|_*.
\end{align}

Because $\left\|\begin{array}{l}
\left[\begin{array}{*{20}{c}}
A &\bm{0}\\
\bm{0} &D
\end{array}\right]
\end{array}\right\|_*
=\|A\|_*+\|D\|_*
=\|\Lambda_A\|_*+\|\Lambda_D\|_*$, to see whether $\left\|\begin{array}{l}
\left[\begin{array}{*{20}{c}}
A &\bm{0}\\
\bm{0} &D
\end{array}\right]
\end{array}\right\|_*=\left\|\begin{array}{l}
\left[\begin{array}{*{20}{c}}
A &B\\
C& D
\end{array}\right]
\end{array}\right\|_*$,
we only need to see whether
\begin{align}\label{Eq:NuclearNormEqualityNewStatement}
\left\|
\begin{array}{l}
\left[\begin{array}{*{20}{c}}
{\begin{array}{l}
\begin{array}{*{20}{c}}
\Lambda_A &\bm{0}\\
\bm{0} &\bm{0}
\end{array}
\end{array}} & B'\\
C' &
{\begin{array}{l}
\begin{array}{*{20}{c}}
\Lambda_D &\bm{0}\\
\bm{0} &\bm{0}
\end{array}
\end{array}}
\end{array}\right]
\end{array}
\right\|_*
=
\left\|
\begin{array}{l}
\left[\begin{array}{*{20}{c}}
{\begin{array}{l}
\begin{array}{*{20}{c}}
\Lambda_A &\bm{0}\\
\bm{0} &\bm{0}
\end{array}
\end{array}} & \bm{0}\\
\bm{0} &
{\begin{array}{l}
\begin{array}{*{20}{c}}
\Lambda_D &\bm{0}\\
\bm{0} &\bm{0}
\end{array}
\end{array}}
\end{array}\right]
\end{array}
\right\|_*
=\|\Lambda_A\|_*+\|\Lambda_D\|_*.
\end{align}

Let the SVD of $E$ be $U_E\Lambda_EV_E^*$, then by Lemma \ref{Lem:NuclearNormEqualityCharacterization-SingularSpace}, the equality holds in (\ref{Eq:NuclearNormEqualityNewStatement}), if and only if there exists a matrix $Q$ such that
\begin{align}
Q^*Q=I,
\end{align}
\begin{align}
U_E=\begin{array}{l}
\left[\begin{array}{*{20}{c}}
I_{r_A\times r_A} &\ &\ &\ \\
\ &\bm{0} &\ &\ \\
\ &\ &I_{r_D\times r_D} &\ \\
\ &\ &\ &\bm{0}
\end{array}\right]
\end{array}
Q,
\end{align}
and
\begin{align}
V_E=
\begin{array}{l}
\left[\begin{array}{*{20}{c}}
I_{r_A\times r_A} &\ &\ &\ \\
\ &\bm{0} &\ &\ \\
\ &\ &I_{r_D\times r_D} &\ \\
\ &\ &\ &\bm{0}
\end{array}\right]
\end{array}
Q,
\end{align}
where $Q$ has only nonzero rows over the indices corresponding to $r_A+r_D$ nonzero columns of
\begin{align}
\begin{array}{l}
\left[\begin{array}{*{20}{c}}
I_{r_A\times r_A} &\ &\ &\ \\
\ &\bm{0} &\ &\ \\
\ &\ &I_{r_D\times r_D} &\ \\
\ &\ &\ &\bm{0}
\end{array}\right]
\end{array}.
\end{align}

Then we can write
\begin{align}
E=\begin{array}{l}
\left[\begin{array}{*{20}{c}}
I_{r_A\times r_A} &\ &\ &\ \\
\ &\bm{0} &\ &\ \\
\ &\ &I_{r_D\times r_D} &\ \\
\ &\ &\ &\bm{0}
\end{array}\right]
\end{array}
Q\Lambda_EQ^*
\begin{array}{l}
\left[\begin{array}{*{20}{c}}
I_{r_A\times r_A} &\ &\ &\ \\
\ &\bm{0} &\ &\ \\
\ &\ &I_{r_D\times r_D} &\ \\
\ &\ &\ &\bm{0}
\end{array}\right]
\end{array}^*.
\end{align}

We can further see that such a matrix $Q$ exists, if and only if the following three conditions hold simultaneously: 1) the row space of $E$ is a subspace of the row space of
\begin{align}
\begin{array}{l}
\left[\begin{array}{*{20}{c}}
I_{r_A\times r_A} &\ &\ &\ \\
\ &\bm{0} &\ &\ \\
\ &\ &I_{r_D\times r_D} &\ \\
\ &\ &\ &\bm{0}
\end{array}\right]
\end{array};
\end{align}
2) the column space of $E$ is a subspace of the column space of
\begin{align}
\begin{array}{l}
\left[\begin{array}{*{20}{c}}
I_{r_A\times r_A} &\ &\ &\ \\
\ &\bm{0} &\ &\ \\
\ &\ &I_{r_D\times r_D} &\ \\
\ &\ &\ &\bm{0}
\end{array}\right]
\end{array};
\end{align}
3) $E$ is a symmetric positive semidefinite matrix.

This proves the claims in this lemma.

\section{Proof for Lemma \ref{Lem:NulcearNormAndTrace}}

Consider the block matrix
\begin{align}
A=
\begin{array}{l}
\left[\begin{array}{*{20}{c}}
A_{1:(n-1),1:(n-1)} &A_{1:(n-1),n}\\
A_{n,1:(n-1)} &A_{n,n}
\end{array}\right]
\end{array}
\end{align}
where $A_{I,J}$ corresponds to submatrix of $A$ with row indices in $I$ and column indices in $J$.

From Lemma \ref{lem:NuclearNormEqualityCharacterization-RowAndColumnSpace}, we have
\begin{align*}
\|A\|_*
& \geq \|A_{1:(n-1),1:(n-1)}\|_*+|A_{n,n}|.
\end{align*}
By applying inductions, we can obtain
\begin{align*}
\|A\|_*\geq \sum_{i=1}^n |A_{i,i}|.
\end{align*}
So $\|A\|_*=\sum_{i=1}^n A_{i,i}$ only if the diagonal elements of $A$ are nonnegative. Moreover, $\|A\|_*=\sum_{i=1}^n A_{i,i}$ if and only if
\begin{align}
\|A\|_*=\left\|\begin{array}{l}
\left[\begin{array}{*{20}{c}}
A_{1,1} &\ &\ &\ \\
\ &A_{2,2} &\ &\ \\
\ &\ &\ddots &\ \\
\ &\ &\ &A_{n,n}
\end{array}\right]
\end{array}\right\|_*.
\end{align}

This means that the nuclear norm of $A$ should be equal to nuclear norm of diagonal matrix
$$\begin{array}{l}
\left[\begin{array}{*{20}{c}}
A_{1,1} &\ &\ &\ \\
\ &A_{2,2} &\ &\ \\
\ &\ &\ddots &\ \\
\ &\ &\ &A_{n,n}
\end{array}\right]
\end{array}.$$

Through similar arguments as in the proof of Theorem \ref{Lem:NuclearNormEqualityCharacterization-SingularSpace} (now extending to multiple blocks), we know that $\|A\|_*=\sum_{i=1}^n A_{i,i}$ if and only if $A=U\Lambda U^*$, where $U$ is a unitary matrix, and $\Lambda$ is a diagonal matrix with nonnegative elements. Then $A$ must be symmetry (or Hermitian in the field of complex numbers) positive semidefinite matrix.

Moreover, when $A$ is a positive semidefinite matrix, we must have $\|A\|_*=\sum_{i=1}^n A_{i,i}$.

\bibliography{Reference_NullSpaceCondition.bib}
\bibliographystyle{IEEEtran}

\end{document}